\documentclass{amsart}
\usepackage{amssymb,amsmath,epsfig}
\usepackage[mathscr]{eucal}

\newtheorem{theorem}{Theorem}[section]
\newtheorem{lemma}[theorem]{Lemma}
\newtheorem{corollary}[theorem]{Corollary}
\newtheorem{proposition}[theorem]{Proposition}
\newtheorem{definition}[theorem]{Definition}
\newtheorem{example}[theorem]{Example}

\numberwithin{equation}{section}
\newcommand{\per}[3]{{\int_{#1}^{#2}\frac{\sin x\cdot #3}{\sin u\sqrt{\sin^2u-\sin^2x}}\,du}}
\newcommand{\ti}[1]{{\hat{#1}}}
\newcommand{\R}{\mathbb{R}}
\newcommand{\Q}{\mathbb{Q}}
\newcommand{\Z}{\mathbb{Z}}
\begin{document}

\title[Geodesics on Orbifolds of Revolution]{On The Existence of Infinitely Many Closed Geodesics on Orbifolds of Revolution}
\thanks{This material is based upon work partially supported by the National Science Foundation under Grant No. DMS-0353622.}
\thanks{The first author wishes to thank M. Stankus for help in producing the figures in this manuscript.}

\author[J. Borzellino]{Joseph E. Borzellino}
\address{Department of Mathematics, California Polytechnic State University, San Luis Obispo, CA 93407 USA}
\email{jborzell@calpoly.edu}

\author[C. Jordan-Squire]{Christopher R. Jordan-Squire}
\address{Department of Mathematics and Statistics, Swarthmore College, Swarthmore, PA 19081 USA}
\email{cjordan1@swarthmore.edu}

\author[G. Petrics]{Gregory C. Petrics}
\address{Department of Mathematics, Middlebury College, Middlebury, VT 05753 USA}
\email{gpetrics@middlebury.edu}

\author[D. Sullivan]{D. Mark Sullivan}
\address{Department of Mathematics, California Institute of Technology, Pasadena, CA 91125 USA}
\email{dsulliva@caltech.edu}


\begin{abstract}
Using the theory of geodesics on surfaces of revolution, we introduce the period function. We use this as our main
tool in showing that any two-dimensional orbifold of revolution homeomorphic to $S^2$ 
must contain an infinite number of geometrically distinct
closed geodesics. Since any such orbifold of revolution can be regarded as a topological two-sphere with metric singularities, we will have extended Bangert's theorem on the existence of infinitely many closed geodesics on any smooth Riemannian two-sphere. In addition, we give an example of a two-sphere cone-manifold of revolution which possesses a single closed geodesic,  thus showing that Bangert's result does not hold in the wider class of closed surfaces with cone manifold structures.
\end{abstract}

\maketitle

\section{Introduction}

In this paper, we study closed geodesics on surfaces of revolution with certain types of metric singularities.
In particular, we are interested in closed (compact, without boundary) surfaces of revolution that are
Riemannian 2-orbifolds. Loosely speaking, an
$2$-orbifold is modeled locally by convex Riemannian surfaces modulo finite groups of
isometries acting with possible fixed points. This means that a neighborhood of each point $p$ of such an orbifold 
is isometric to a Riemannian quotient $U_p/\Gamma_p$ where $U_p$ is a convex Riemannian surface diffeomorphic to $\mathbb{R}^2$,
and $\Gamma_p$ is a finite group of isometries acting effectively on $U_p$.
Every Riemannian surface is trivially an orbifold, with each $\Gamma_p$
being the trivial group.  The reader interested in more background on orbifolds
should consult \cite{BorPHD}, Thurston's classic \cite{Thurston78}, or the more recent textbook \cite{MR1299730}.
For the purposes of this paper, however, we will only need to apply a simple explicit criterion
to determine whether a closed surface of revolution is a 2-orbifold (see
section~\ref{OrbsOfRevSection}).

The existence of closed geodesics on Riemannian manifolds has a long and storied past dating back to Poincar\'{e} 
\cite{MR1729907}. It seems that not much has been done on the existence of closed geodesics in singular spaces. The existence of at least one closed geodesic on a compact 2-orbifold was shown in \cite{MR1419471} and closed geodesics in orbifolds of higher dimensions have recently been studied in \cite{GuruprasadHaefliger}. The paper \cite{NorburyRubinstein}
studies the issue of closed geodesics in spaces with incomplete metrics. The relevance here is that a complete Riemannian orbifold with singular set removed is a Riemannian manifold with incomplete metric and it is known \cite{MR1218706} that closed geodesics
in a complete Riemannian orbifold may not pass through the singular set, unless they are entirely contained within it.

Here we are interested in the question of the existence of infinitely many closed geodesics.
In \cite{MR1209957}, Bangert used the
work of Franks \cite{MR1161099} to show that every smooth Riemannian $S^2$
has infinitely many closed geodesics. For orbifolds with $S^2$ as
the underlying topological space, the existence of an infinity of
closed geodesics is an open question. In the general category of closed surfaces
of revolution with singular points (which have underlying
topological space $S^2$), one may construct examples with exactly one closed
geodesic (see example~\ref{void}), showing that analogue of Bangert's result is false in this category.
We call such a surface {\em void}.
A spherical 2-orbifold of revolution is a closed two-dimensional surface of revolution homeomorphic to $S^2$ that satisfies
a certain special orbifold condition at its north and south poles. It is natural to ask whether
void orbifolds of revolution exist. In resolving this question we extend Bangert's result by proving that

\begin{theorem}\label{MainTheorem} Every spherical 2-orbifold of
revolution has infinitely many closed geodesics.
\end{theorem}

Since we are dealing only with surfaces of revolution, our
techniques are relatively elementary. We begin by recalling the basic theory
about surfaces of revolution and their geodesics, most of which can
be found in \cite{MR0394451}, \cite{MR0203595}, or \cite{OpreaJ}. 

\section{Basic Theory}\label{BasicTheory}

In what follows the term {\em smooth function} will refer to a function of class $C^\infty$. 
In fact, $C^2$ is sufficient for our needs.

\begin{definition}\label{SurfaceDef} Let $\alpha:[u_N,u_S]\to\R^2$ be a simple (no self intersections) smooth
plane curve $\alpha(u)=(g(u),h(u))$  where $g$ and $h$ are smooth functions on the interval $[u_N,u_S]$,
with $h\geq 0$, and $h(u)=0$ if and only if $u=u_N$ or $u=u_S$.  A {\em spherical surface of revolution} $M$ is a
surface embedded
isometrically in $\R^3$ that admits a parametrization
$\mathbf{x}:[u_N,u_S]\times\R\to M$ of the form
$$\mathbf{x}(u,v)=(g(u),h(u)\cos v,h(u)\sin v),$$
That is, $M$ is the surface of revolution obtained by rotating $\alpha$ about the $x$-axis.
The curve $\alpha$ will be called the {\em profile curve}.
\end{definition}

Note that a spherical surface of revolution $M$ is necessarily homeomorphic to $S^2$ and that by definition the sets 
$N=\mathbf{x}(u_N,v)$ and $S=\mathbf{x}(u_S,v)$ for $v\in\R$ reduce to single points
which will be referred to as the {\em north} and {\em south poles} of
$M$. Metric singularities may only occur at these two points. $M$ is smooth everywhere else. We also do not require that
the function $g$ be monotone.
Throughout this paper all surfaces of revolution will be assumed spherical as in definition~\ref{SurfaceDef} even though much of the classical theory we review applies equally well to any surface of revolution.

Rotation about the $x$-axis in $\R^3$ descends to a natural $S^1$-action $S^1\times M\to M$ on $M$ by isometries:
\[
(e^{i\theta},(x,y,z))\mapsto
(x,y\cos\theta-z\sin\theta,y\sin\theta+z\cos\theta).
\]
This action is free except at the north and
south poles which remain fixed.

For a surface of revolution $M$, a simple computation gives the
coefficients of the first fundamental form or metric tensor (subscripts denote partial derivatives):
$$
E=\mathbf{x}_u\cdot\mathbf{x}_u=[g'(u)]^2+[h'(u)]^2,\quad
F=\mathbf{x}_u\cdot\mathbf{x}_v=0 \quad\textrm{and}\quad
G=\mathbf{x}_v\cdot\mathbf{x}_v=h^2(u),
$$
so that the metric (away from any singular point) is
$$ds^2=\left([g'(u)]^2+[h'(u)]^2\right)du^2+h^2(u)dv^2.$$  

Note that the parametrization is orthogonal ($F=0$) and that $E_v=G_v=0$. 
Surfaces given by parametrizations with these properties are said to be {\em $u$-Clairaut}.

For any $u$-Clairaut surface, and hence any surface of
revolution, the geodesic equations reduce to
\begin{equation}
\label{geq1} u'' + \frac{E_u}{2E}{u'}^2-\frac{G_u}{2E}
{v'}^2=0
\end{equation}
\begin{equation}
\label{geq2} v'' + \frac{G_u}{G} u' v' = 0.
\end{equation}

A curve $\gamma(t)=\mathbf{x}(u(t),v(t))$ on $M$ is a geodesic if and only
if the above equations are satisfied by the coordinate functions $u$ and $v$
of $\gamma$.  Also, a geodesic satisfying these
equations must be parametrized proportional to arc length and hence has constant
speed. In particular, we may assume that $\gamma$ has unit speed. That is,
$\|\gamma'\|=Eu'^2+Gv'^2\equiv 1$. The existence and uniqueness theorem for solutions
of ordinary differential equations implies that, given a
point in $p$ in $M$ and a vector $\xi$ in $T_pM$, the tangent
plane to $M$ at $p$, there is a unique geodesic $\gamma$ satisfying
$\gamma(0)=p$ and $\gamma'(0)=\xi$.

We now recall two important classes of geodesics on surfaces of
revolution.

\begin{example} {\em A unit speed curve $\gamma(t)=\mathbf{x}(u(t),v(t))$ with
$v(t)\equiv v_0$, a constant, is a $u$-parameter curve or {\em meridional arc}. Such curves are always geodesics.
To see this, note that $v'=v''\equiv 0$, so equation (\ref{geq2}) is satisfied trivially.  The
unit speed relation is, in this case, $E{u'}^2=1$, so ${u'}^2=1/E$.
Differentiating each side and dividing by $2u'$ gives
\[
u'' = -\frac{E_uu'+E_vv'}{2u' E^2} = -\frac{E_u}{2E^2}
=-\frac{E_u}{2E}{u'}^2,
\]
which is equivalent to (\ref{geq1}), since $v'\equiv 0$.  We will use the term
{\em meridian} for those meridional arcs that join $N$ to $S$. }
\end{example}

\begin{example} {\em A unit speed curve $\gamma(t)=\mathbf{x}(u(t),v(t))$ with
$u(t)\equiv u_0\in (u_N,u_S)$, a constant, is a $v$-parameter curve or {\em parallel arc}.
For a parallel arc we have $u'=u''\equiv 0$ and the unit speed relation $G{v'}^2=1$.
Differentiating the unit speed relation yields $v''=-\frac12G_v/G^2=0$ since $G=h(u)^2$ depends only on $u$.
Thus (\ref{geq2}) is satisfied. Equation (\ref{geq1}) reduces to
$G_u{v'}^2=0$.
Hence parallel arcs are geodesic precisely when $G_u(u_0)=0$, or
equivalently, $h'(u_0)=0$. We will use the term {\em parallel} for those parallel arcs which are entire circles. }
\end{example}

For the remainder of the paper we will assume all geodesics come with unit speed parametrizations.

The main classical tool used to get qualitative information about geodesics on surfaces of revolution is the
{\em Clairaut relation}, which we present now. Let $\gamma(t)$ be a geodesic on $M$. Then
$\gamma'=u'\mathbf{x}_u+v'\mathbf{x}_v$. If there exists $t_0$ with $v'(t_0)=0$, then the uniqueness of geodesics implies that
$\gamma$ must be a meridional arc as $\gamma'$ is parallel to $\mathbf{x}_u$ at $t_0$. As a result $v'$ cannot change sign, and we may assume, without loss of generality that $v'(t)\ge 0$. In fact, $v'(t)>0$ unless $\gamma$ is a meridional arc.

Let $\varphi_\gamma(t)=\angle(\gamma',\mathbf{x}_u)$ be the angle between
$\gamma'$ and $\mathbf{x}_u$ at time $t$. Since the surface parametrization $\mathbf{x}$ and $\gamma$ are smooth, 
$\varphi_\gamma(t)$ is a smooth function that takes its values in the interval $[0,\pi]$. From the discussion above
we see that $\varphi_\gamma(t)\in (0,\pi)$ for all $t$ if and only if $\gamma$ is not a meridional arc.
Now consider the quantity $c_{\gamma}=Gv'$ along a geodesic
$\gamma$. Then
$$c_{\gamma}'=Gv''+(G_uu' +G_vv')v'=Gv''+G_uu' v'=0$$
where the second equality follows since $G_v=0$ and the last equality follows from the
second geodesic equation (\ref{geq2}).
Thus the quantity $c_\gamma$ is constant along geodesic paths.  Comparing the two expressions for 
$\gamma'\cdot\textbf{x}_v$:
\begin{align*}
\gamma'\cdot\textbf{x}_v & =
(u'\textbf{x}_u+v'\textbf{x}_v)\cdot\textbf{x}_v = Gv' = c_{\gamma}\\
\gamma'\cdot\textbf{x}_v  & = \|\gamma'\|\;\|\textbf{x}_v\|
\cos\left(\frac{\pi}{2}-\varphi_{\gamma}\right) =
\sqrt{G}\sin\varphi_{\gamma}
\end{align*}

yields the {\em Clairaut relation}:

\begin{proposition}\label{ClairautProp}
If $c_{\gamma}(t)=G(t)v'(t)$ along a geodesic $\gamma(t)$, then the quantity
\begin{equation}
\label{CR}
c_\gamma(t)=\sqrt{G(t)}\sin\varphi_{\gamma}(t)=h(u(t))\sin\varphi_{\gamma}(t)
\end{equation}
is constant.
\end{proposition}

The constant $c_\gamma$ is called the {\em slant} of $\gamma$.
Since $0\leq\sin\varphi_\gamma(t)\leq 1$ for all $t$
we must have that $h(u(t))\geq c_\gamma$ for all $t$. That is, $\gamma$ is
must lie entirely in the region of the surface $M$ where $h(u)\geq c_\gamma$. 

\begin{corollary}\label{MeridianPoles} For a spherical surface of revolution, a geodesic $\gamma$ with an endpoint at either
pole must be a meridional arc.
\end{corollary}
\begin{proof} Since $\gamma$ has an endpoint at a pole assume for concreteness that $\gamma(a)=N$ and that $\gamma$ is defined over an interval $[a,b]$.
Let $t_n\to a$ be a sequence of real numbers $t_n\in (a,b)$ converging to $a$. Then $h(u(t_n))\to h(u(a))=0$, whence
$c_\gamma(t_n)\to 0$. By proposition~\ref{ClairautProp}, $c_\gamma(t)\equiv 0$, which implies that
$\sin\varphi_{\gamma}(t)\equiv 0$. Thus, $\gamma$ is a meridional arc.
\end{proof}

\begin{corollary}\label{NonMeridionalExtension} If $\gamma$ is not a meridional arc, then $\gamma$ cannot pass through a pole of $M$. Thus, non-meridional geodesics $\gamma$ have a unique extension to a unit speed geodesic $\hat\gamma:\mathbb{R}\to M$.
\end{corollary}
\begin{proof} This follows from corollary~\ref{MeridianPoles} and the existence and uniqueness theorem for geodesics.
\end{proof}

\section{Qualitative Theory and a Classification of Geodesics}

In light of corollary~\ref{NonMeridionalExtension} we will now assume that all non-meridional geodesics will be defined
on $\mathbb{R}$, and meridional arcs are extended to meridians.

Motivated by the Clairaut relation we define, for $c>0$, the super-level sets 
$$M^c=\{p\in M\mid\text{if }p=\mathbf{x}(u,v),\text{then } h(u)>c\}$$

Points of $M^c$ will be referred to as points of $M$ with $h(u)>c$ for convenience.
$M^c$ may have several connected components $M^c_{(j)}$, but it is always true that the boundary of any such
component $\partial M^c_{(j)}=\rho_0\cup\rho_1$ where $\rho_i$ are parallels 
$\rho_i: t\mapsto\mathbf{x}(u_i,t)$ with $h(u_i)=c$ for $i=0,1$. 
Note also that if $c$ is a local minimum value of $h(u)$ then $c>0$ and the closure
$\overline{M^c_{(j)}}$ is a proper subset of $\{p\in M\mid\text{if }p=\mathbf{x}(u,v),\text{then } h(u)\ge c\}$.

The qualitative behavior of non-meridional geodesics given next is key to our analysis.

\begin{proposition}\label{OscillateAsymptoticGeodesics} Suppose a geodesic $\gamma(t)=\mathbf{x}(u(t),v(t))$
on a spherical surface of revolution $M$ is tangent at $t=t_0$ 
to a non-geodesic parallel $\rho_0:t\mapsto\mathbf{x}(u_0,t)$ of $M$.  Then $\gamma$ is constrained to lie in the connected component $\overline{M_\gamma}$ of $\overline{M^{c_\gamma}}$ which contains $\gamma$. The boundary 
$\partial\overline{M_\gamma}=\rho_0\cup\rho_1$ where 
$\rho_1:t\mapsto\mathbf{x}(u_1,t)$ is a parallel of $M$ with $h(u_1)=c_\gamma$. 
Moreover, $\gamma$ either oscillates between the parallels $\rho_i$ intersecting them tangentially or $\gamma$ spirals asymptotically to $\rho_1$ which is necessarily a geodesic parallel.
\end{proposition}
\begin{proof} If $\gamma$ is tangent at $t=t_0$ 
to a non-geodesic parallel $\rho_0$, then $u(t_0)=u_0$ and $\gamma'(t_0)$ is parallel to $\mathbf{x}_v$. Thus,
$\varphi_\gamma(t_0)=\pi/2$ which implies $c_\gamma=h(u(t_0))=h(u_0)$.
Thus, by the Clairaut relation we may then conclude that the entire geodesic $\gamma$ 
lies in a region of $M$ that corresponds to points where the profile curve is $\ge c_\gamma$. Since $\gamma$ is not a parallel, (otherwise, $\gamma$ would have to coincide with $\rho_0$ which is not geodesic), $h(u(t))>c_\gamma$ for some $t\in\R$. Thus, 
$\gamma$ is a subset of a connected component $\overline{M_\gamma}$ of $\overline{M^{c_\gamma}}$. Since $\rho_0$ is non-geodesic,
$u_0$ is not a critical point of $h$ and thus $h$ is monotone in a neighborhood of $u=u_0$. The Clairaut relation then implies that
$\gamma$ lies on one side of $\rho_0$. That is, for all $t$, $u(t)\in (u_N,u_0]$ or $u(t)\in [u_0,u_S)$. This shows that
$\rho_0\subset\partial\overline{M_\gamma}$ and that intersections of $\gamma$ with $\rho_0$ are tangential. 

Without loss of generality, we may assume that $u(t)\in [u_0,u_S)$. Let $u_1\in (u_0,u_S)$ be the smallest number such that 
$h(u_1)=c_\gamma$. Let $\rho_1:t\mapsto\mathbf{x}(u_1,t)$ be the corresponding parallel.
Define $b=\sup_{t\in\R}u(t)$.  

If there is a $t_b$ such that $u(t_b)=b$, then $u'(t_b)=0$. 
Thus, $\gamma$ is parallel to $\mathbf{x}_v$ at $t_b$ and hence 
$c_\gamma=h(b)$. By the choice of $u_1$, we must have $b=u_1$ and thus,
$\gamma\cap\rho_1\ne\emptyset$. As before, we may conclude that $\rho_1$ is non-geodesic, all intersections are tangential, 
and $\gamma$ lies on one side of $\rho_1$. 
In particular, the set $\{u(t)\mid t\in\R\}=[u_0,u_1]$ and $\gamma$ oscillates back and forth between the two parallels $\rho_0$ and $\rho_1$ which form
the boundary $\partial\overline{M_\gamma}$.

On the other hand, if no such $t_b$ exists, then $\lim_{t\to\infty}u(t)=b$ and $\gamma$ is asymptotic to the parallel at $u=b$. Since $\gamma$ is geodesic, this implies that the parallel $\rho_b$ at $u=b$ is geodesic with slant $c_{\rho_b}=h(b)$. By taking a limits we
conclude that  $c_\gamma=\lim_{t\to\infty}\big[h(u(t))\sin\angle(\gamma'(t),\mathbf{x}_u)\big]=h(b)\sin\angle(\rho_b',\mathbf{x}_u)=h(b)$. 
By the choice of $u_1$ we conclude that $b=u_1$ and that $\rho_b=\rho_1$. In particular, in this case, $\gamma$ spirals asymptotically to a geodesic parallel $\rho_1$ and $\partial\overline{M_\gamma}=\rho_0\cup\rho_1$.
\end{proof} 

Geodesics which exhibit the oscillating behavior of proposition~\ref{OscillateAsymptoticGeodesics} will be called {\em oscillating geodesics} and those with asymptotic behavior will be called {\em asymptotic geodesics}. There is actually one last type of geodesic, called a {\em bi-asymptotic geodesic}. This is
a geodesic that spirals into a geodesic parallel as $t\to -\infty$ and another geodesic parallel as $t\to\infty$. The existence of bi-asymptotic geodesics will be considered in proposition~\ref{existasymp} where we consider conditions that imply the existence of (bi)-asymptotic geodesics.

\begin{proposition}\label{existasymp}  Let $\alpha=(g,h):[u_N,u_S]\to\R^2$ be the profile curve of $M$. 
Then $\Gamma_M$ contains an asymptotic geodesic if and only if $h$ has a critical point in the
interval $(u_N,u_S)$ that is not a local maximum.
\end{proposition}

\begin{proof}  Suppose $h$ has a critical point at $u=u_0$ that is not a local maximum. Without loss of generality, we may assume there is
$u_1\in h^{-1}(h(u_0))$ such that the open interval $(u_0,u_1)\subset\{u\in(u_N,u_S)\mid h(u)>h(u_0)\}$. 

If $h'(u_1)\ne 0$ then the corresponding parallel
at $u_1$ is non-geodesic and by proposition~\ref{OscillateAsymptoticGeodesics} there is a geodesic $\gamma$ through $\mathbf{x}(u_1,0)$ and parallel to
$\mathbf{x}_v(u_1,0)$ so that $\gamma$ is asymptotic to the parallel at $u_0$.

On the other hand, if $h'(u_1)=0$, then pick a point $\hat{u}\in (u_0,u_1)$. Since $h(u_0)<h(\hat{u})$ we can find $\hat{\varphi}\in (0,\pi/2)$ with
$\hat{\varphi}=\displaystyle\arcsin\left(\frac{h(u_0)}{h(\hat{u})}\right)$. Now, let $\gamma$ be the geodesic with 
$\gamma(0)=\mathbf{x}(\hat{u},0)$ and with
$\angle(\gamma'(0),\mathbf{x}_u)=\hat{\varphi}$. Then the slant of $\gamma$, $c_\gamma=h(\hat{u})\sin\hat{\varphi}=h(u_0)$. By proposition~\ref{OscillateAsymptoticGeodesics}, we may conclude that $\gamma$ is asymptotic to
the geodesic parallels at $u_0$ and $u_1$. In this case, $\gamma$ is bi-asymptotic.
\end{proof}

We have shown that a geodesic on spherical surface of revolution is either a meridian, a geodesic parallel, an oscillating geodesic, an asymptotic geodesic or a bi-asymptotic geodesic. We now define the boundary values and boundary function on the set of
geodesics on $M$.

\begin{definition} Let $\gamma(t)=\mathbf{x}(u(t),v(t))$ be  a non-meridional geodesic
on a spherical surface of revolution $M$. Define $b_0(\gamma)=\inf_{t\in\R}(u(t))$ and 
$b_1(\gamma)=\sup_{t\in\R}(u(t))$ to be the {\em left} and 
{\em right} {\em boundary values} of $\gamma$, respectively. 
If $\gamma$ is a meridian we define $b_0(\gamma)=u_N$ and 
$b_1(\gamma)=u_S$. 
The {\em boundary function} $b$ of $\gamma$ is defined by
$b(\gamma)=(b_0(\gamma),b_1(\gamma))\in [u_N,u_S]\times [u_N,u_S]$.
\end{definition}

In the case of a non-meridional geodesic, by proposition~\ref{OscillateAsymptoticGeodesics} we have that
$b_0(\gamma)$ and $b_1(\gamma)$ are the corresponding $u$ values for the parallels $\rho_0$ and $\rho_1$.
When the geodesic under consideration is clear, we will often drop the reference to $\gamma$ and
refer to the boundary values of $\gamma$ as $b_0$ and $b_1$.

\begin{definition}\label{equivalence}
Two geodesics $\gamma_1$ and $\gamma_2$ on $M$ are
{\em equivalent} if $\gamma_1$ and $\gamma_2$ are in the same orbit of the
natural $S^1$ action on $M$.
We denote the set of all equivalence classes $[\gamma]$ by $\Gamma_M$.
\end{definition}

Since the $S^1$ action on $M$ preserves parallels, we conclude from proposition~\ref{OscillateAsymptoticGeodesics} that
the boundary function $b:\Gamma_M\to [u_N,u_S]\times [u_N,u_S]$ is well-defined and injective.
We adopt the common abuse of notation by simply referring to a
geodesic $\gamma\in\Gamma_M$.

We can classify all geodesics on a spherical surface of revolution by boundary function:

\begin{definition}\label{classification}
Let $\gamma$ be a geodesic in $\Gamma_M$.
\begin{enumerate}
\item If $b=(u_N,u_S)$ then $\gamma$ is a meridian. Otherwise, if $\gamma$ is not a meridian,
\item If $b_0=b_1$, $\gamma$ is a geodesic parallel,
\item If $b_0 \neq b_1$ with $h'(b_0)$ and $h'(b_1)$ both non-zero, then $\gamma$ is an oscillating geodesic.
\item If $b_0 \neq b_1$  with both $h'(b_0)$ and $h'(b_1)$ equal to zero, then $\gamma$ is a bi-asymptotic geodesic, otherwise
\item If $b_0 \neq b_1$  with either $h'(b_0)$ or $h'(b_1)$ equal to zero, then $\gamma$ is an asymptotic geodesic.
\end{enumerate}
\end{definition}

\begin{figure}[ht]
   \centering
   \scalebox{0.3}{
        \epsfig{file=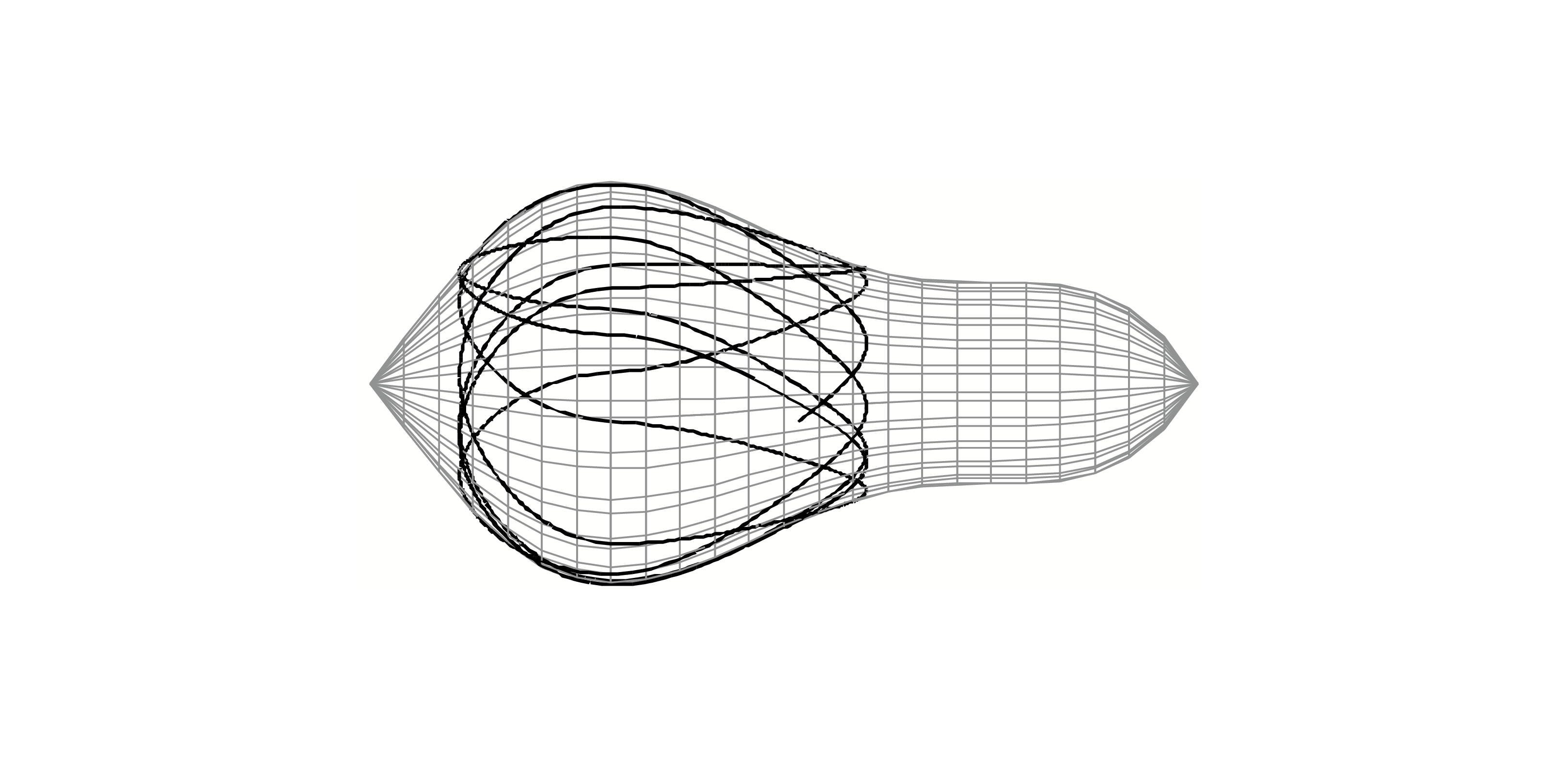}}
   \scalebox{0.24}{
        \epsfig{file=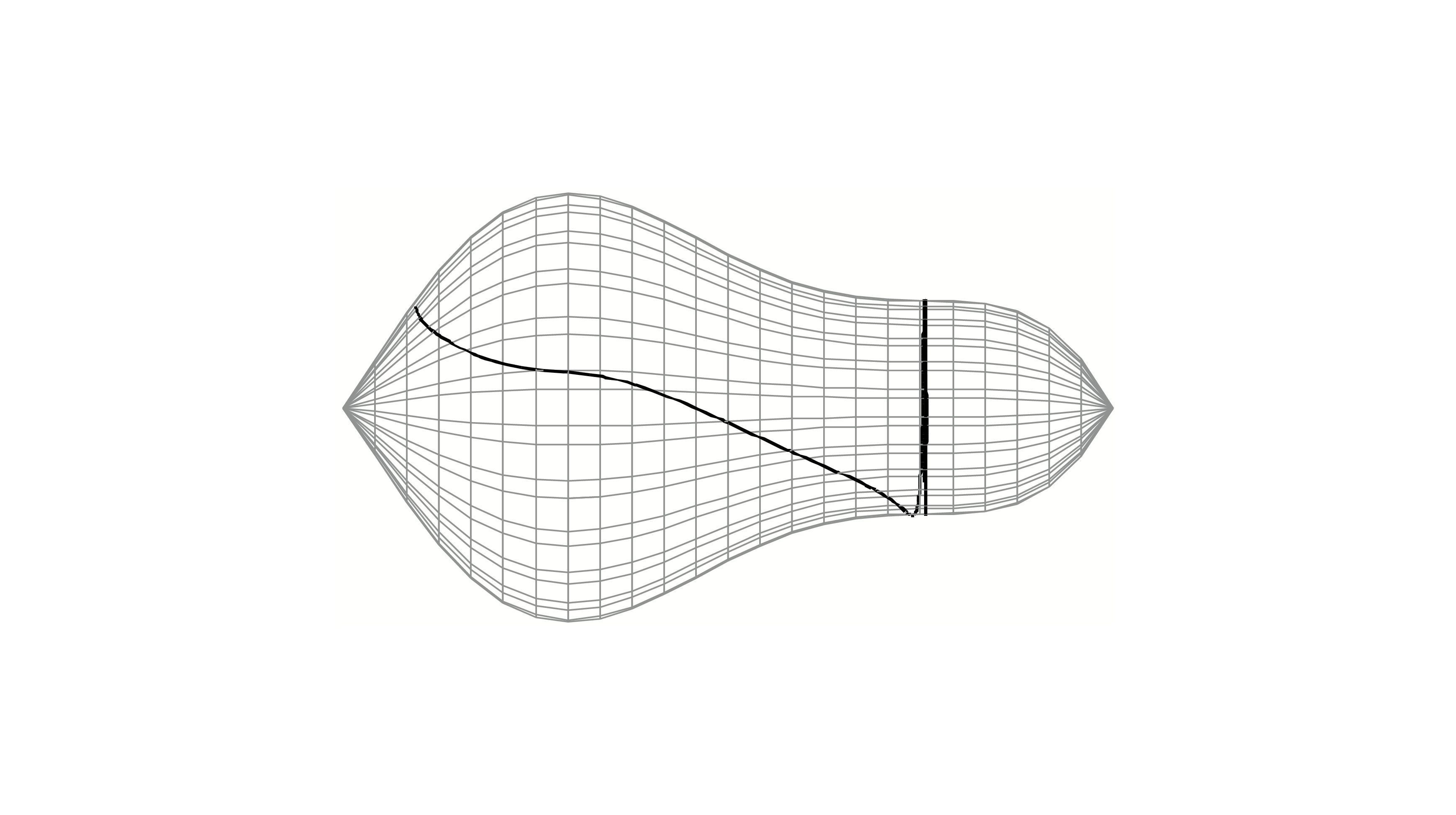}}
\parbox{3in}{\caption{\label{GeodesicFigure} An oscillating and asymptotic geodesic on a spherical surface of revolution}}
\end{figure}

Examples of an oscillating and asymptotic geodesics are given in figure~\ref{GeodesicFigure}.

Let $\gamma$ be an oscillating geodesic. If $u(t_0)=b_0=u(t_1)$ for 
$t_0\neq t_1$ and there is a unique $t\in (t_0,t_1)$ such that $u(t)=b_1$, we
call the segment of $\gamma$ corresponding to the interval
$[t_0,t_1]$ an {\em oscillation}. Since it is not important for what follows, we will refer to 
bi-asymptotic geodesics as simply asymptotic geodesics also.

\section{A Topology on the Set of Oscillating Geodesics}

\begin{definition} A geodesic $\gamma$ is {\em closed} if there
exist real numbers $t_0\neq t_1$ such that
$$\gamma(t_0)=\gamma(t_1)\quad\textup{and}\quad\gamma'(t_0)=
\gamma'(t_1).$$
\end{definition}

Equality of the derivatives distinguish closed geodesics from the more general notion of geodesic loop.
Every geodesic parallel is closed, and no asymptotic geodesic or meridian (using our definition) is closed.
Oscillating geodesics, however, may or may not be closed. Since we are interested in closed geodesics
the set
$$\Gamma_M^O = \{ [\gamma] : \gamma \textup{ is an oscillating
geodesic}\}$$ will be the most interesting to us.

Note that if $\gamma$ is oscillating, then $\gamma$ is the unique
geodesic with left boundary $b_0(\gamma)$. This is because
$h'(b_0(\gamma))\ne 0$, so the parallel at $u=b_0(\gamma)$ is not geodesic
and there can be no geodesic asymptotic to it. Thus, by our classification,
any geodesic which shares a left boundary with $\gamma$ must be oscillating itself.
But, any oscillating geodesic intersects its left boundary tangentially, so by the definition
of our equivalence relation and the uniqueness of geodesics we conclude that $\gamma$ is the unique
geodesic in its equivalence class with left boundary $b_0(\gamma)$.
Thus, the map $b_0:\Gamma_M^O\to (u_N,u_S)$ is injective.  In
particular, for oscillating geodesics, the right boundary value is
determined by the left boundary value.

\begin{proposition}\label{TopologyonGammaO} Let $b_1(u_1)=\min\{u>u_1:h(u)=h(u_1)\}$ and let 
$\mathscr{U}=\{u_1\in(u_N,u_S):h'(u_1)>0 \textup{ and } h'(b_1(u_1))<0\}$.  Then
$\mathscr{U}$ is an open subset of the interval $(u_N,u_S)$ and $b_0:\Gamma_M^O\to\mathscr{U}$
is a bijection.
\end{proposition}

\begin{proof}  We first show that $b_0$ is a bijection. Indeed, $b_0(\gamma)\in \mathscr{U}$ for any
$\gamma\in\Gamma_M^O$ by proposition~\ref{OscillateAsymptoticGeodesics}.  
For any $u_1\in \mathscr{U}$, there is a
geodesic $\gamma$ with the initial conditions $u(0)=u_1$,
$u'(0)=0$.  Then $h'(u_1)>0$ implies $\gamma$ is not a geodesic parallel and $b_0(\gamma)=u_1$.
Thus, $b_1(\gamma)=b_1(u_1)$, and $h'(b_1(u_1))<0$ implies $\gamma$ is not asymptotic, so $\gamma\in\Gamma_M^O$.

Smoothness of $h$ implies that if $u_1\in \mathscr{U}$ there exists an $\epsilon>0$
such that $(u_1-\epsilon,u_1+\epsilon)\subset \mathscr{U}$. Thus,
$\mathscr{U}$ is an open subset of the real interval $(u_N,u_S)$.
\end{proof}

We can now topologize on $\Gamma_M^O$ by pull-back:
We declare a subset $U\subset\Gamma_M^O$ to be open if and only if $b_0(U)$ is open
in $\mathscr{U}$. Hence $\Gamma_M^O$ is
homeomorphic to a disjoint union of open intervals of $(u_N,u_S)$. 
This allows us, for example, to speak
of a sequence of geodesics in the space $\Gamma_M^O$ as a sequence
of (left) boundary values from $\mathscr{U}$.  Consequently, we can also
easily define convergence of oscillating geodesics, and, more
importantly, continuous functions defined on $\Gamma_M^O$.

\section{The Period Function}\label{PeriodFunctionSection}

We now present our main analytic tool for detecting closed geodesics on
spherical surfaces of revolution.

In the case of oscillating or asymptotic geodesics, the geodesic equations
(\ref{geq1}) and (\ref{geq2}) can be reduced to a first-order system
and solved explicitly. Equation (\ref{geq2}), after dividing by
$v'$ (which is never zero for oscillating or asymptotic geodesics)
and integrating, becomes
$$
\int\frac{v''}{v'}\;dt = -\int\frac{G_u}{G}u'\;dt
    =-\int\frac{G'}{G}\;dt,
$$
since $G_v=0$.  Then $v'=c/G$ for some constant $c\in\R$
and $c=Gv'=c_\gamma$, again showing the slant $c_\gamma$ to be constant.
Using $v'=c_\gamma/G$ in the unit speed relation
$Eu'^2+Gv'^2=1$ gives
$$
u'=\pm\sqrt{\frac{G-c_\gamma^2}{EG}}.
$$
Hence,
$$
\frac{dv}{du}=\frac{v'}{u'}=\pm\frac{c_\gamma\sqrt{E}}{\sqrt{G}
{\sqrt{G-c_\gamma^2}}},
$$
and
$$
v=\pm\int\frac{c_\gamma\sqrt{E}}{\sqrt{G}{\sqrt{G-c_\gamma^2}}}\;du.
$$
As we will soon see, by measuring the total change in $v$ an oscillating geodesic makes
between its boundaries one can determine if it is closed.
This motivates the following definition.

\begin{definition}
The {\em period function} $\Phi_{M} : \Gamma_M^O \rightarrow (0,\infty)$ is defined by
$$
\Phi_{M}(\gamma)=2\int_{b_0(\gamma)}^{b_1(\gamma)}\frac{c_\gamma\sqrt{E}}{\sqrt{G}\sqrt{G-c_\gamma^2}}\;du =%
2\int_{b_0(\gamma)}^{b_1(\gamma)}\frac{c_\gamma\sqrt{E}}{h(u)\sqrt{h^2(u)-c_\gamma^2}}\;du.
$$
We denote the integrand by $f_\gamma(u)$.
\end{definition}

Geometrically, the period function gives the change in $v$ as
$\gamma$ undergoes one oscillation.  Since
$h^2(b_0)=h^2(b_1)=c_{\gamma}^2$, the integral is improper for every
geodesic $\gamma$, however, because it represents the change in $v$
between $b_0$ and $b_1$ it must converge for every
$\gamma\in\Gamma^O_M$.  We can use this geometric
interpretation to see that the period function is invariant under reparametrization and scaling of $M$ and
to extend the domain of the period function to
include the asymptotic geodesics, by setting
$\Phi_M(\gamma_0)=\infty$ for any asymptotic geodesic $\gamma_0$.

The next theorem shows how the period function can be used to detect closed geodesics.

\begin{theorem}\label{closed}
An oscillating geodesic $\gamma$ on a spherical surface of revolution $M$ is
closed if and only if $\Phi_M(\gamma)=2q\pi$ for some rational
number $q\in\mathbb{Q}$.
\end{theorem}

\begin{proof}
We can assume without loss of generality that $\gamma$ satisfies the initial conditions
$$
\gamma(0) = \mathbf{x}(b_0,0) \quad\mbox{and}\quad \gamma'(0)=\frac{\mathbf{x}_v(b_0,0)}{\|\mathbf{x}_v(b_0,0)\|}.
$$
If $\gamma(t)=\mathbf{x}(u(t),v(t))$ is closed, there exists a
$t_0\in\R^+$ such that $\gamma(t_0)=\gamma(0)$ and
$\gamma'(t_0)=\gamma'(0)$.  In particular,
$\gamma(t_0)=\mathbf{x}(b_0,2r\pi)$ for some positive integer $r$.
Note that the period function does not depend on the value $v(0)$,
so by rotational symmetry, $v$ changes the same amount during every
oscillation of $\gamma$.  Clearly, between $t=0$ and $t=t_0$,  $\gamma$ has completed, say, $s$ oscillations. That is, 
there have been $s$ times subsequent to $t=0$ that $u(t)$ has re-attained the boundary value $b_0$.  
Therefore, $\Phi_M(\gamma)=2(r/s)\pi$.

Conversely, suppose $\Phi_M(\gamma)=2(r/s)\pi$ for some
$r,s\in\Z^+$, where $\gamma$ is taken to have the same initial
conditions.  Then there exists a $t_0\in\R^+$ such that
$$
\gamma(t_0)=\mathbf{x}(b_0,2r\pi)=\mathbf{x}(b_0,0)=\gamma(0).
$$
Since $u(t_0)=b_0$, by proposition~\ref{OscillateAsymptoticGeodesics}, we must have $\gamma'(t_0)$ tangent to $\mathbf{x}_v(b_0,0)$, and thus,
$\gamma'(t_0)=\gamma'(0)$. Hence $\gamma$ is closed.
\end{proof}

The next theorem shows that the period function is continuous.  
A sketch of the proof of this result first appeared in the unpublished manuscript \cite{REU2004}.
For clarity of the exposition, we relegate to section~\ref{PeriodProof} the rather technical proof of this result.

\begin{theorem}\label{cont}
If $\gamma_0\in\Gamma_M^O$, then $\Phi_M$ is continuous at $\gamma_0$.
\end{theorem}

As we will see, the continuity of $\Phi_M$ at every oscillating
geodesic implies the existence of infinitely many geodesics on many
spherical surfaces of revolution.

\begin{definition}
A non-empty open subset $U\subseteq\Gamma_M^O$ on which $\Phi_M$ is
a constant, irrational multiple of $\pi$ is said to be {\em void}.
\end{definition}

The definition is motivated by theorem~\ref{closed}, which implies that all
oscillating geodesics $\gamma$ with $b_0(\gamma)\in U$ fail to close smoothly, so $U$
is {\em void} of closed geodesics.

\begin{corollary}\label{nonvoid}
Suppose a spherical surface of revolution $M$ has a non-empty open subset $U$ of $\Gamma_M^O$ that is
not void. Then $M$ has infinitely many closed geodesics.
\end{corollary}

\begin{proof}
Let $\Phi_M(U)=\{\Phi_M(\gamma):\gamma\in U\}$.
If $\Phi_M(U)$ is a constant rational multiple of $\pi$ we are done by theorem~\ref{closed}, so suppose $\Phi_M$
is not constant on $U$.  By continuity of $\Phi_M$, there exists a
nonempty open interval $I\subset\Phi_M(O)$.  $\Q\pi$ is
dense in any such $I$ yielding an infinite number of closed geodesics in $U$ by theorem~\ref{closed}.
\end{proof}

The following corollary shows that the existence of an
asymptotic geodesic on $M$ implies the existence of such a non-void subset
of $\Gamma^O_M$, and hence the existence of infinitely many closed
geodesics.

\begin{corollary}\label{AsympImpliesClosed} Let $M$ be a spherical surface of revolution
with an asymptotic geodesic $\gamma_0$ asymptotic to the geodesic parallel
at $b_0(\gamma_0)$. Then if $\gamma_n\to\gamma_0$ is a sequence of oscillating geodesics,
$$
\lim_{n\rightarrow\infty}\Phi_M(\gamma_n)=\Phi_M(\gamma_0)=\infty.
$$
Thus, by corollary~\ref{nonvoid}, $M$ has infinitely many closed geodesics.
\end{corollary}

\begin{proof}
Let $A>0$. $\Phi_M(\gamma_0)=\displaystyle\int_{b_0(\gamma_0)}^{b_1(\gamma_0)}f_{\gamma_0}=\infty$, so there exists $\delta,\mu>0$ so that
$\displaystyle A<\int_{b_0(\gamma_0)+\delta}^{b_1(\gamma_0)-\mu}f_{\gamma_0}$. Choose $N>0$ large enough so that $b_0(\gamma_n)<b_0(\gamma_0)+\delta$ and
$b_1(\gamma_n)>b_1(\gamma_0)-\mu$ for $n>N$. Thus, 
$$
\Phi_M(\gamma_n)=\int_{b_0(\gamma_n)}^{b_1(\gamma_n)}f_{\gamma_n}>\int_{b_0(\gamma_0)+\delta}^{b_1(\gamma_0)-\mu}f_{\gamma_n}
$$
On the interval $(b_0(\gamma_0)+\delta,b_1(\gamma_0)-\mu)$, $f_{\gamma_n}\to f_{\gamma_0}$, and both $f_{\gamma_n}$ and $f_{\gamma_0}$ are
bounded hence integrable. Thus, by dominated convergence, for $\varepsilon>0$, there is $N'>N$ so that
$$
\Phi_M(\gamma_n)>\int_{b_0(\gamma_0)+\delta}^{b_1(\gamma_0)-\mu}f_{\gamma_n}>\left[\int_{b_0(\gamma_0)+\delta}^{b_1(\gamma_0)-\mu}f_{\gamma_0}\right]%
-\varepsilon>A-\varepsilon
$$
for $n>N'$. This implies $\Phi_M(\gamma_n)\to\infty$.
\end{proof}

\begin{corollary}\label{MultipleCPsImplyClosed} A spherical surface of revolution whose profile curve has more than one
critical point necessarily has an infinite number of closed geodesics.
\end{corollary}
\begin{proof} This follows from proposition~\ref{existasymp} and corollary~\ref{AsympImpliesClosed}.
\end{proof}

\begin{corollary}\label{SingleCPImplyOneOrMany} A spherical surface of revolution whose profile curve has a single critical point (which must be
a maximum), has exactly one closed geodesic or infinitely many.
\end{corollary}
\begin{proof} The parallel at the critical point is necessarily geodesic. If $\Gamma^O_M$ is not void, then there
are an infinite number of closed geodesics by corollary~\ref{nonvoid}. The only other possibility is that $\Phi_M$ is a constant irrational multiple of $\pi$
over its entire domain $\Gamma^O_M$. Then by theorem~\ref{closed}, no oscillating geodesic is closed, and $M$ has exactly one closed geodesic.
\end{proof}

\begin{definition}\label{VoidSurface} A spherical surface of revolution with exactly one closed geodesic will be called a {\em void surface}.
\end{definition}
An explicit example of a void surface will be given in section~\ref{TwoExamplesSection}.

\section{Surfaces of Revolution with Constant Period Function}

Since, ultimately, we wish to show that no void spherical 2-orbifolds of revolution exist, corollary \ref{SingleCPImplyOneOrMany} implies we should look
for general conditions that imply the period function is constant. We do exactly that in this section.

If a spherical surface of revolution $\textbf{x}(u,v)=(g(u),h(u)\cos v,h(u)\sin v)$ obtained from the profile curve $\alpha(u)=(g(u),h(u))$ is to have a constant period function, we can, without loss of generality, assume that 
$h(u)$ is a smooth function from $[0,L]$ to $[0,1]$
satisfying:
\begin{enumerate}
\item $h(0)=h(L)=0$
\item $h$ has a unique critical point, say $u_0$, on $[0,L]$
\item $h(u_0)=1$
\end{enumerate}
where
$$g(u)=\int_0^u\sqrt{1-{[h'(t)]}^2}\,dt.$$

As a result, we may assume the metric on $M$ is of the
form
$$
ds^2=du^2+h^2(u)dv^2.
$$

If the period function $\Phi_M$ is to be constant, proposition \ref{existasymp} and corollary~\ref{AsympImpliesClosed} imply that condition (2)
is necessary. (1) and (3) may be
satisfied by an appropriate reparametrization and scaling of the profile curve, which
does not affect the period function. The following proposition is adapted from \cite{MR496885}.

\begin{proposition}\label{MetricHasSpecialForm}
Let $M$ be a spherical surface of revolution
satisfying conditions (1),(2) and (3).  We can define new
coordinates $(r,v)$ on $M$ so that the metric in these coordinates has the form
$$
ds^2=E(r)dr^2+\sin^2r\;dv^2,
$$
where $\ti{E}(\cos r)=E(r)$ is a function from $[0,\pi]$ to $\mathbb{R}^+$.
\end{proposition}

\begin{proof}
Define $f:[0,L]\rightarrow [0,\pi]$ and $c:[-1,1]\rightarrow [0,L]$
by
$$
f(u)=\begin{cases}
 \arcsin h(u) & \text{if $u\in[0,u_0]$}\\
 \pi-\arcsin h(u) & \text{if $u\in[u_0,L]$}
\end{cases}
$$
$$
c(\cos r)=\begin{cases}
 (h|_{[0,u_0]})^{-1}(\sin r) & \text{if $r\in [0,\frac{\pi}{2}]$}\\
 (h|_{[u_0,L]})^{-1}(\sin r) & \text{if $r\in [\frac{\pi}{2},\pi]$.}
\end{cases}
$$

Then, setting $r=f(u)$, we have $c(\cos r)=f^{-1}(r)=f^{-1}(f(u))=u$
and
\[
h(u)=h[c(\cos r)]=\sin r.
\]
Hence,
\[
du^2=\left(\frac{1}{f'(u)}\right)^2dr^2=\left(\frac{\cos r}{h'[c(\cos r)]}\right)^2dr^2,
\]
and we can now write the metric on $M$ as
\[
ds^2=E(r)dr^2+\sin^2r\;dv^2,
\]
where $\ti{E}(\cos r)=E(r)$ is a function from $[0,\pi]$ to $\mathbb{R}^+$
defined by
$$E(r)=
\begin{cases}
 \ti{E}(\cos r)=\displaystyle\frac{\cos^2 r}{\big(h'[c(\cos r)]\big)^2} & \text{if $r\neq\frac{\pi}{2}$}\\
 \\[-.15in]
 \ti{E}(0)=\displaystyle\frac{1}{{[f'(u_0)]}^2}=\frac{-1}{h''(u_0)} & \text{if $r=\frac{\pi}{2}$}.
\end{cases}
$$
Note that condition (2) implies $h''(u_0)<0$ and that by differentiating the relation $h(u)=\sin f(u)$ twice and
evaluating at $u=u_0$ shows that $f'(u_0)=\sqrt{-h''(u_0)}$. Thus $E$ is continuous on $[0,\pi]$.
\end{proof}

We can now take as a starting point in our search for surfaces with
constant constant period function those surfaces of revolution with
metric of the form
$$
ds^2=E(u)\,du^2+\sin^2 u\,dv^2,
$$
where $E(u)$ is a function from $[0,\pi]$ to $\R^+$.  This
corresponds to the spherical surface of revolution $M$ with profile curve
$\alpha(u)=(g(u),\sin u)$, where
$$
g(u)=\int_0^u\sqrt{E(t)-\cos^2t}\,dt.
$$
If $\gamma_x$ is the geodesic with left boundary value $b_0=x$, then
the right boundary value $b_1=\pi-x$ and the period function may then be written as a
function of $x\in(0,\frac{\pi}{2})$:
$$
\Phi_M(\gamma_x)=\Phi_M(x)=\per{x}{\pi-x}{\sqrt{E}(u)},
$$
which is continuous on $(0,\frac{\pi}{2})$ by theorem~\ref{cont}.  The following technical
lemma adapted from \cite{MR496885} will be essential in our characterization of surfaces of
revolution with constant period function.

\begin{lemma}\label{TechnicalLemma} Consider the function
$$
F(x)=\per{x}{\pi-x}{f(u)}
$$
Define a function $\ti{f}$ by the formula $f(u)=\ti{f}(\cos u)$. Then $F(x)$ is identically zero on $(0,\frac{\pi}{2})$ if and only if
$\ti{f}$ is an odd function over $[-1,1]$.
\end{lemma}

\begin{proof}
Let $\ti{f}^e(\cos u)=(\ti{f}(\cos u)+\ti{f}(-\cos u))/2$ be the
even part of $\ti{f}$.  Then $\ti{f}$ is odd if and only if
$\ti{f}^e$ is identically zero.  We have
\begin{equation*}
\begin{split}
F(x) &= \per{x}{\frac{\pi}{2}}{\ti{f}(\cos u)}+\per{\frac{\pi}{2}}{\pi-x}{\ti{f}(\cos u)}\\
     &= \per{x}{\frac{\pi}{2}}{\ti{f}(\cos u)}\\
     & \quad\quad +\per{\frac{\pi}{2}}{\pi-x}{2\ti{f}^e(\cos u)}-\per{\frac{\pi}{2}}{\pi-x}{\ti{f}(-\cos u)}\\
     &= \per{x}{\frac{\pi}{2}}{\ti{f}(\cos u)}-\per{x}{\frac{\pi}{2}}{\ti{f}(\cos u)}\\
     & \quad\quad +\per{x}{\frac{\pi}{2}}{2\ti{f}^e(\cos u)}\\
     &= 2\per{x}{\frac{\pi}{2}}{\ti{f}^e(\cos u)}.
\end{split}
\end{equation*}

So $\ti{f}^e(\cos u)\equiv 0$ gives $F(x)=0$ for all
$x\in (0,\frac{\pi}{2})$.\\

For the converse, we follow a proof given in \cite{MR496885}.  Assume
$F(x)$ is zero for all $x\in(0,\frac{\pi}{2})$, then
$$
F(x)=\per{x}{\pi-x}{f(u)}=\per{x}{\frac{\pi}{2}}{2\ti{f}^e(\cos u)}\equiv 0.
$$
So the function
$$
I(a)=\int_a^{\frac{\pi}{2}}\frac{\cos x\cdot F(x)}{\sqrt{\sin^2x-\sin^2a}}dx
$$
is zero for $a\in(0,\frac{\pi}{2}]$.  Also, for such $a$, the function
$$
\frac{1}{\left(\sqrt{\sin^2u-\sin^2x}\right)\left(\sqrt{\sin^2x-\sin^2a}\right)}
$$
is (Lebesgue) integrable on the set
$\{(u,x)\in [a,\frac{\pi}{2}]\times [a,\frac{\pi}{2}]\mid x\leq u\}$.
Applying Fubini's Theorem we have
$$
I(a) = \int_a^{\frac{\pi}{2}}\frac{2\ti{f}^e(\cos u)}{\sin u}\left(\int_a^u\frac{\sin x\cos x}{\sqrt{\sin^2u-\sin^2x}%
\sqrt{\sin^2x-\sin^2a}}dx\right)du.
$$
The substitution $t=\sqrt{\sin^2x-\sin^2a}\big/\sqrt{\sin^2u-\sin^2x}$ gives
\begin{equation*}
\begin{split}
I(a) &= 2\int_a^{\frac{\pi}{2}}\frac{\ti{f}^e(\cos u)}{\sin u}\left(\int_0^{\infty}\frac{dt}{1+t^2}\right)du=%
\pi\int_a^{\frac{\pi}{2}}\frac{\ti{f}^e(\cos u)}{\sin u}du.
\end{split}
\end{equation*}
As $\sin u$ is strictly positive on $(0,\frac{\pi}{2}]$, $I(a)=0$
for all $a\in[0,\frac{\pi}{2}]$ implies $\ti{f}^e(\cos u)=0$ for all
$\cos u\in[-1,1]$. That is, $\ti{f}$ is odd on $[-1,1]$.
\end{proof}

\begin{proposition}\label{TechnicalProposition} For a spherical surface of revolution M with metric
$$
ds^2=E(u)\,du^2+\sin^2(u)\,dv^2,
$$
define $a_c(u)=\sqrt{E}(u)-c$ for any $c\in\R^+$. Then
$\Phi_M(x)\equiv 2c\pi$ on $(0,\frac{\pi}{2})$ if and only if
the function $\ti{a}_c$ defined by $\ti{a}_c(\cos u)=a_c(u)$ is an odd function from $[-1,1]$ to
$[-c,c]$.
\end{proposition}

\begin{proof} Let $S^2$ be the standard 2-sphere of constant curvature 1 in $\R^3$
generated as a surface of revolution by the profile curve
$\alpha(u)=(\cos u,\sin u)$. Thus, $S^2$ is parametrized by
$$
\textbf{x}(u,v)=(\cos u,\sin u\cos v,\sin u\sin v).
$$
The geodesics on $S^2$ are great circles, so
$$
\Phi_{S^2}(x)=2{{\int_{x}^{\pi-x}\frac{\sin x}{\sin u\sqrt{\sin^2u-\sin^2x}}\;du}}\equiv 2\pi.
$$

Then, for all $x\in (0,\frac{\pi}{2})$,
\begin{equation*}
\begin{split}
\Phi_M(x)&=2\per{x}{\pi-x}{\sqrt{E}(u)}
    = 2\per{x}{\pi-x}{(c+\ti{a}_c(\cos u))}\\
    &= c\Phi_{S^2}(x)+2\per{x}{\pi-x}{\ti{a}_c(\cos
    u)}\\
    &= 2c\pi+2\per{x}{\pi-x}{\ti{a}_c(\cos u)}.
\end{split}
\end{equation*}
The proof of the proposition now follows from lemma~\ref{TechnicalLemma}, which implies that
$\ti{a}_c$ must be odd.  For $u\in (0,\pi)$, $c+\ti{a}_c(\cos u)=\sqrt{E}(u)>0$ so $\ti{a}_c(\cos u)>-c$ for $u\in (0,\pi)$. This implies that
$\ti{a}_c(-\cos u)>-c$, so since $\ti{a}_c$ is odd, we have $\ti{a}_c(\cos u)=a_c(u)\in [-c,c]$ for $u\in (0,\pi)$.
\end{proof}

At this point we are able to recover Bangert's result for spherical surfaces of revolution which have (smooth) Riemannian metrics, such as
ellipsoids of revolution. We first need a computation.

Let $\phi_N$, resp. $\phi_S$, be the angle between the profile curve
$\alpha(u)=(g(u),h(u))=(g(u),\sin u)$ and the axis of rotation at $g(0)$, resp.
$g(\pi)$.  Then
\begin{subequations}\label{Phi_NSEq}
\begin{gather}
\label{Phi_NEq}\sin\phi_N=
\frac{h'(0)}{\sqrt{{[g'(0)]}^2+{[h'(0)]}^2}}= \frac{\cos(0)}{\sqrt{E(0)}}=\frac{1}{c+\ti{a}_c(1)}\\
\intertext{and}
\label{Phi_SEq}\sin\phi_S=
\frac{-h'(\pi)}{\sqrt{{[g'(\pi)]}^2+[{h'(\pi)]}^2}}=
\frac{-\cos(\pi)}{\sqrt{E(\pi)}}=
\frac{1}{c+\ti{a}_c(-1)}=\frac{1}{c-\ti{a}_c(1)},
\end{gather}
\end{subequations}
with the last equality following since $\ti{a}_c$ is odd on $[-1,1]$.

We now easily deduce Bangert's result for Riemannian spherical surfaces of revolution.

\begin{corollary}\label{SpecialCaseBangert}
Every smooth Riemannian $S^2$ arising as a surface of revolution has infinitely
many closed geodesics.
\end{corollary}
\begin{proof}
The result follows if the surface has non-constant period function by corollary~\ref{nonvoid}.
Thus, we assume the surface has constant period function. Since the surface is a smooth manifold, the
profile curve meets the $x$-axis at right angles, so that
$\sin\phi_N=\sin\phi_S=1$. Equations~\eqref{Phi_NEq} and \eqref{Phi_SEq} imply that $c+\ti{a}_c(1)=c-\ti{a}_c(1)=1$ so
$0=\ti{a}_c(1)=\ti{a}_c(-1)$ and $c=1$.  Hence $\Phi_M\equiv 2\pi$
and all oscillating geodesics close up after one oscillation.
\end{proof}

\section{Orbifolds of Revolution}\label{OrbsOfRevSection}

Our work up to this point is valid for spherical surfaces of revolution in general. Since our main theorem~\ref{MainTheorem} concerns orbifolds,
we now specialize to that case. Spherical orbifolds of revolution are easily identifiable by their tangent cones at the poles. Namely, the tangent cone
at a pole must be isometric to the metric quotient of the flat plane $\R^2$ by a finite cyclic group of rotations fixing the origin. 
Note that the tangent cone at a pole is generated by rotating the tangent line to the profile curve at the pole about the axis of rotation. If the cyclic groups at the poles are of different orders, the orbifold is commonly referred to as {\em bad} since it will not arise as a quotient of a Riemannian $S^2$ by a finite cyclic group of isometries \cite{Thurston78}.

\begin{figure}[ht]
   \centering
   \scalebox{0.4}{
        \epsfig{file=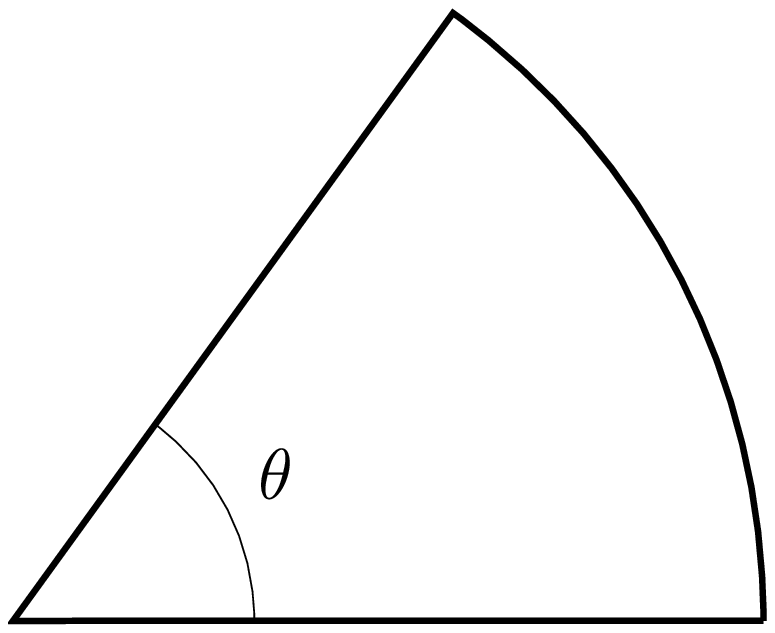}}\raisebox{.6in}{$\qquad\leadsto\qquad$}
        \scalebox{0.32}{
        \epsfig{file=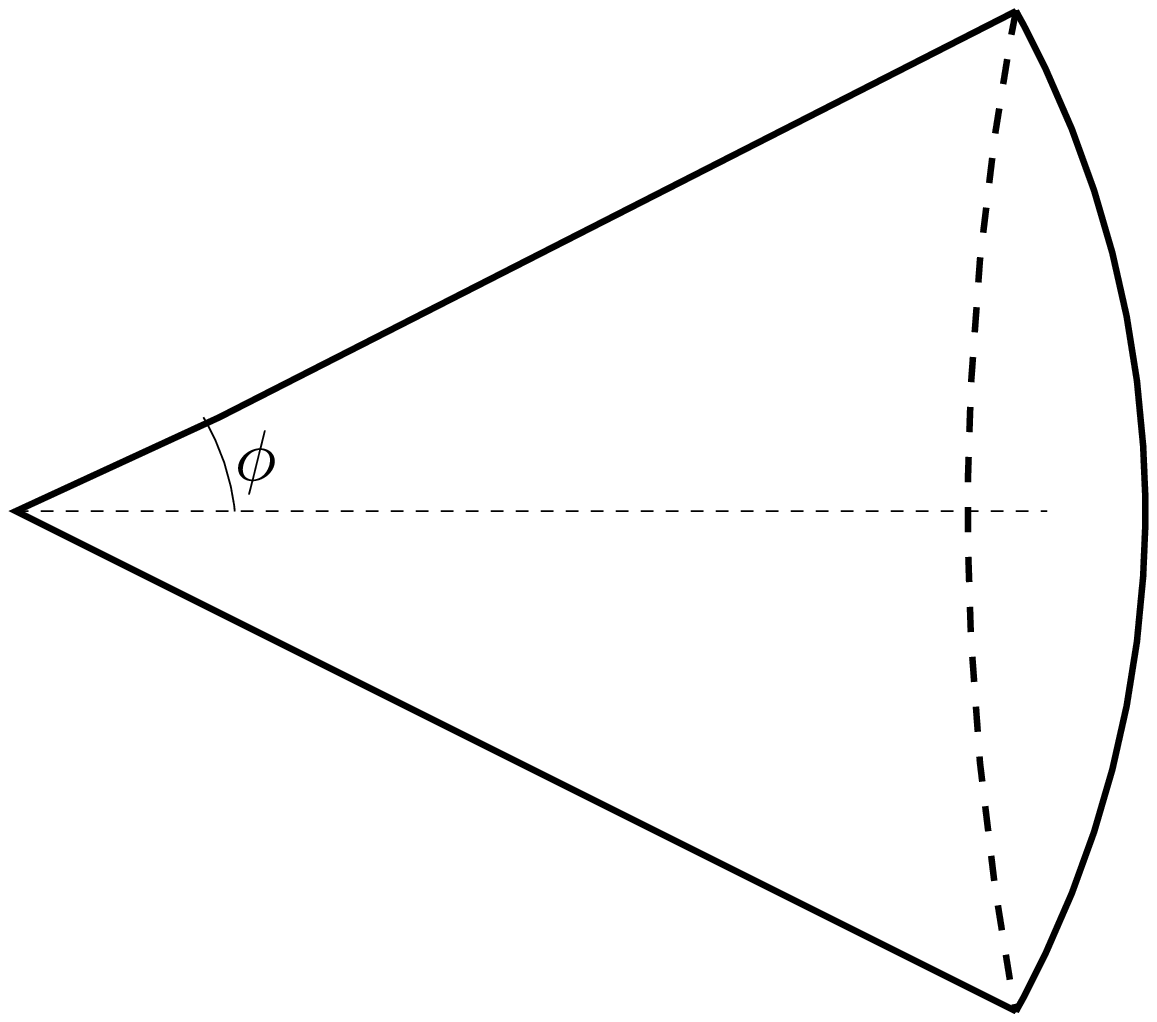}} 
   \caption{\label{WedgeConeFigure} Cone as quotient of a planar sector}
   \end{figure}

In general, a flat right circular cone with vertex angle $\phi$ is obtained by identifying the edges of a plane circular sector of angle $\theta$.
The relation between $\theta$ and $\phi$ is easily computed: $\theta=2\pi\sin\phi$. See figure~\ref{WedgeConeFigure}.
Thus, if the tangent cone at a pole of spherical orbifold of revolution
is isometric to $\R^2/\Z_m$, then $\theta=2\pi/m$ for a positive integer $m$. So, for an orbifold of revolution, if $\phi_N$ and $\phi_S$ are as in equations~\eqref{Phi_NSEq}, we must have $\sin\phi_N = 1/m$ and $\sin\phi_S = 1/k$ for some positive integers $m$ and $k$.

We have the following restriction for spherical orbifolds of revolution of constant period function.

\begin{theorem}\label{crational}
Let M be a spherical orbifold of revolution with metric
$ds^2=E(u)\,du^2+\sin^2(u)\,dv^2$. Then $\Phi_M(x)\equiv 2c\pi$ on
$(0,\frac{\pi}{2})$ implies $c$ is rational.
\end{theorem}

\begin{proof}  Equations~\eqref{Phi_NSEq} give
$$
\sin\phi_N=\frac{1}{c+\ti{a}_c(1)} \quad\mbox{and}\quad
\sin\phi_S=\frac{1}{c+\ti{a}_c(-1)}=\frac{1}{c-\ti{a}_c(1)},
$$
since $\ti{a}_c$ is odd on $[-1,1]$.  As noted above, if $M$ is an orbifold,
then $c+\ti{a}_c(1)$ and $c-\ti{a}_c(1)$ must be integers. This easily implies that $c$ be rational.  
In fact, $c=n/2$ for some
positive integer $n$ and $\Phi_M(x)=n\pi$ on $(0,\frac{\pi}{2})$.
\end{proof}

We are now in a position to prove the main result of this paper.

\begin{theorem}There are no void spherical orbifolds of revolution.
Hence, every orbifold of revolution has infinitely many closed
geodesics.
\end{theorem}

\begin{proof} Suppose one such example existed. By corollary~\ref{MultipleCPsImplyClosed}, we may assume that the profile curve has
a single critical point and hence by proposition~\ref{MetricHasSpecialForm} that the metric on $M$ is of the form required in 
theorem~\ref{crational}. By theorem~\ref{closed} and corollary~\ref{nonvoid},
$\Phi_M$ must be a constant, irrational multiple of $\pi$.  However, by Theorem \ref{crational}, an orbifold
of revolution with constant $\Phi_M$ must have $\Phi_M\equiv 2c\pi$
with $c\in\Q$. Hence no such void spherical orbifold exists and all spherical orbifolds of revolution must have infinitely many geodesics.
\end{proof}

\section{Two examples}\label{TwoExamplesSection}

In summary, we can characterize all spherical surfaces of revolution with constant
period function as having a metric of the form $ds^2=(c+f(\cos u))^2\,du^2+\sin^2(u)\,dv^2$ where
\begin{enumerate}
\item $c$ is a real constant,
\item $f(\cos u)$ is an odd function from $[-1,1]$ to $[-c,c]$.
\end{enumerate}
The void spherical surfaces of revolution satisfy these conditions but have
$c\notin\Q$, and hence are not orbifolds.  The orbifolds of
revolution with constant period function must satisfy (1), (2) and
\begin{enumerate}
\item[(3)] $c+f(1)$ and $c-f(1)$ are positive integers.
\end{enumerate}

\begin{figure}[ht]
   \centering
   \scalebox{0.3}{
        \epsfig{file=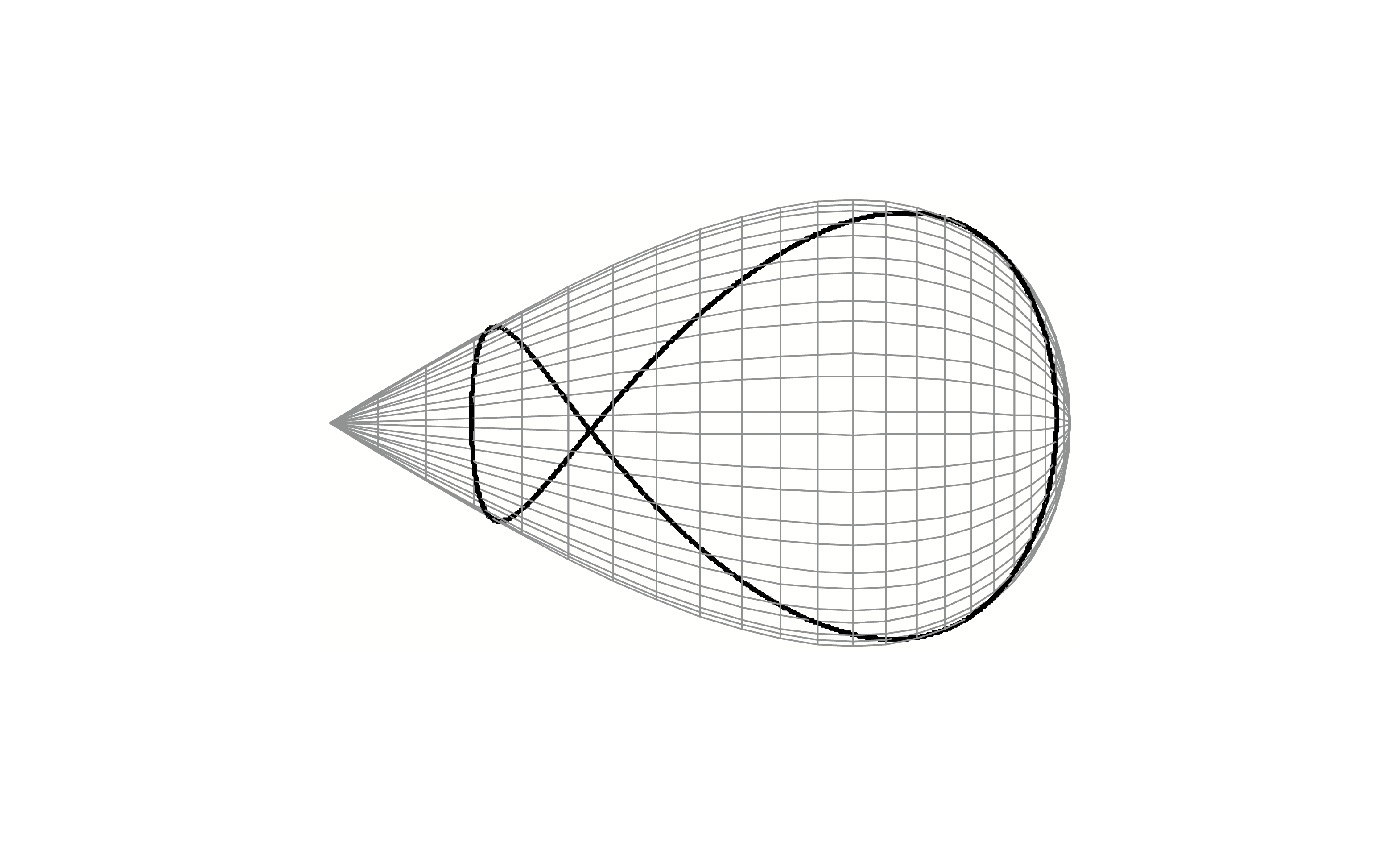}}
   \caption{\label{TanneryPearFigure} A typical closed geodesic on a Tannery pear}
   \end{figure}

\begin{example}[Tannery's pear]
{\em Take $c=2$ and $a_c(u)=\cos u$ (so $\ti{a}_c(\cos u)$ is the
identity map on $[-1,1]$, and hence odd).  This surface, known as
Tannery's pear, has a period function that is constant $4\pi$, so
all non-meridional geodesics are closed.  It also is an orbifold.

Taking as a profile curve $\alpha(u)=(g(u),h(u))$, where $h(u)=\sin u$ and
\begin{equation*}
\begin{split}
g(u) &= \int_0^u\sqrt{E(t)-(h'(t))^2}\;dt= \int_0^u\sqrt{(2+\cos t)^2-\cos^2t}\;dt\\
    &= \int_0^u\sqrt{4+4\cos t}\;dt = 4\sqrt{2}\sin(u/2)
\end{split}
\end{equation*}
gives a parametrization for Tannery's pear in $\R^3$:
$$
\mathbf{x}(u,v) = (4\sqrt{2}\sin(u/2), \sin u\cos v, \sin u\sin v),
$$
where $u\in[0,\pi]$ and $v\in[0,2\pi]$.

We also have that
\[
\sin\phi_N=\frac{\cos(0)}{\sqrt{E(0)}}=\frac{1}{3}\quad\textrm{and}
\quad\sin\phi_S=\frac{-\cos(\pi)}{\sqrt{E(\pi)}}=1,
\]
so Tannery's pear is an orbifold.  In orbifold
terminology, Tannery's pear is a $\Z_3$-teardrop, as the metric
is actually smooth at $u=\pi$ and the single cone point at $u=0$ is
of order 3. See figure~\ref{TanneryPearFigure}.}
\end{example}

\begin{figure}[ht]
   \centering
   \scalebox{0.3}{
        \epsfig{file=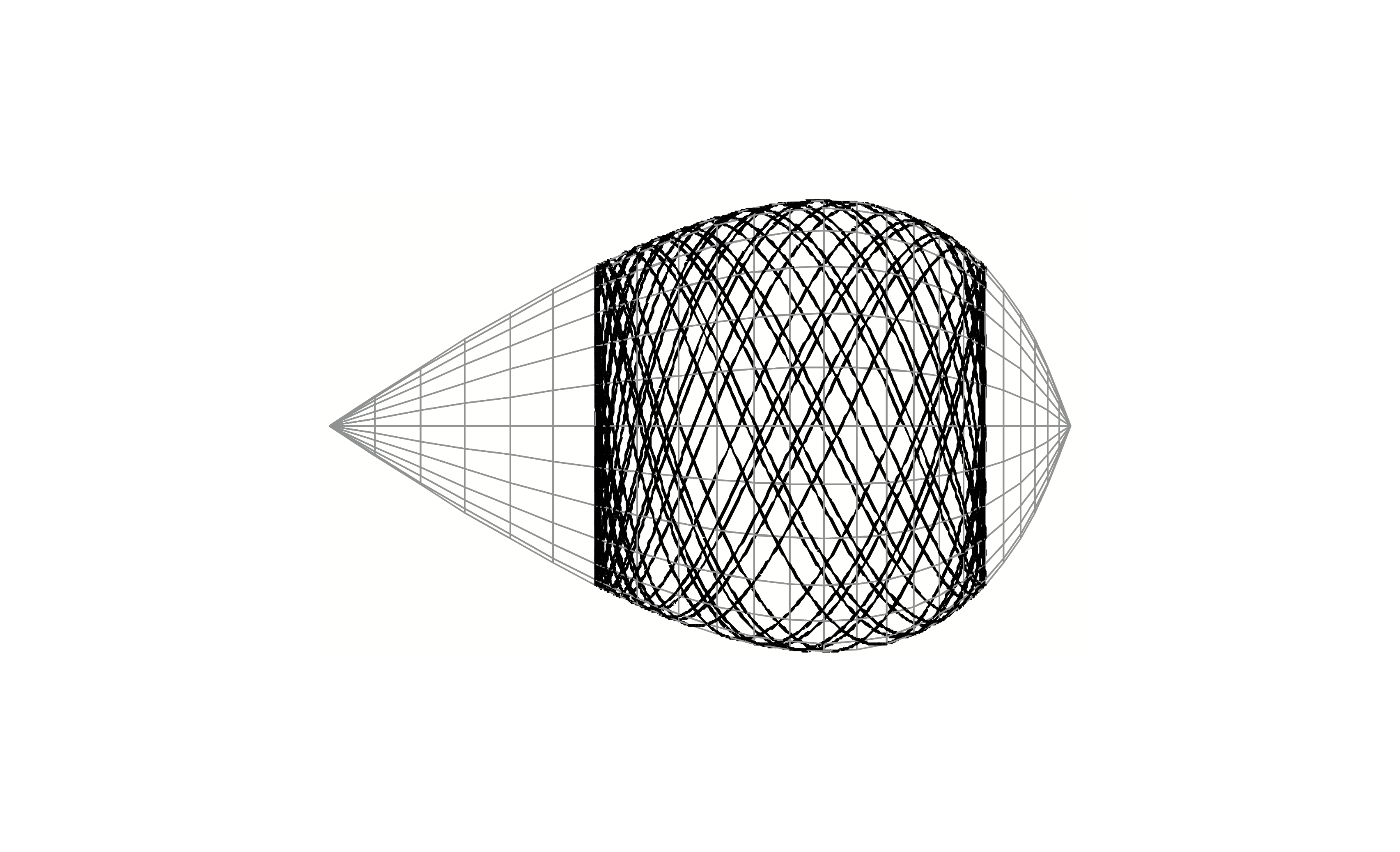}} 
   \caption{\label{VoidSurfaceFigure} An oscillating geodesic that is not closed}
   \end{figure}

\begin{example}[A void surface]\label{void}
{\em In the previous example, take $c=\sqrt{5}$.  A surface with this
metric can be isometrically embedded in $\R^3$ by the
parametrization
$$
\mathbf{x}(u,v)=\left(\int_0^u\sqrt{5+2\sqrt{5}\cos t}\;dt,\sin
u\cos v,\sin u\sin v\right).
$$
It has constant period function $2\pi\sqrt{5}$ and hence, its only
closed geodesic is the parallel at $u=\pi/2$.  However,
$$
\sin\phi_N=\frac{1}{\sqrt{5}+1}\quad\text{and}\quad
\sin\phi_S=\frac{1}{\sqrt{5}-1},
$$
so, like all void spherical surfaces of revolution, this one is not an orbifold.
See figure~\ref{VoidSurfaceFigure}.} 
\end{example}

\section{Proof of the Continuity of the Period Function}\label{PeriodProof}

In this section, we prove theorem~\ref{cont}, which asserts the continuity of the period function. The notation used will be
that from section~\ref{PeriodFunctionSection}.

\begin{proof}
Without loss of generality, we will assume the profile curve $\alpha=(g,h)$ of $M$ is parametrized by arclength. Thus,
$E=[g'(u)]^2+[h'(u)]^2\equiv 1$. Let $\gamma_0$ be an oscillating geodesic. By proposition~\ref{TopologyonGammaO} there is
$0<\varepsilon_0<1$ so that $h'(u)>0$ on $B(b_0(\gamma_0),\varepsilon_0)$ and $h'(u)<0$ on
$B(b_1(\gamma_0),\varepsilon_0)$. Here $B(p,r)$ is the open interval of radius $r$ centered at $p$. By shrinking 
$\varepsilon_0$ if necessary, we may assume $|h'(u)|\ge\eta>0$ on 
$B=B(b_0(\gamma_0),\varepsilon_0)\cup B(b_1(\gamma_0),\varepsilon_0)$.

Let $\gamma_n$ be a sequence of oscillating geodesics with $b_0(\gamma_n)\to b_0(\gamma_0)$ and
$b_0(\gamma_n)\in B(b_0(\gamma_0),\frac14\varepsilon_0^2)$. We may also assume
$b_1(\gamma_n)\in B(b_1(\gamma_0),\frac14\varepsilon_0^2)$, by choosing $n$ large enough.

Consider the integrand $f_{\gamma_n}(u)$ of the period function $\Phi_M$:
$$
f_{\gamma_n}(u) = \frac{c_{\gamma_n}}{h(u)\sqrt{h^2(u)-c_{\gamma_n}^2}}
=\frac{c_{\gamma_n}}{h(u)\sqrt{h(u)+c_{\gamma_n}}}\cdot\frac{1}{\sqrt{h(u)-c_{\gamma_n}}}
$$
for $u\in(b_0(\gamma_n),b_1(\gamma_n))$. Since $c_{\gamma_n}=h(b_0(\gamma_n))$, applying the mean value theorem to the second factor in the last equality yields:
$$
f_{\gamma_n}(u)=\frac{h(b_0(\gamma_n))}{h(u)\sqrt{h(u)+h(b_0(\gamma_n))}}\cdot\frac{1}{\sqrt{h'(\xi_u)}}%
\frac{1}{\sqrt{u-b_0(\gamma_n)}}
$$
for some $\xi_u\in (b_0(\gamma_n),u)$. Let $m=\inf_B h$ and define $\lambda=(2m\eta)^{-1/2}$. Then since
$h(b_0(\gamma_n))<h(u)$ for $u\in (b_0(\gamma_n),b_1(\gamma_n))$, we have
\begin{align*}
f_{\gamma_n}(u) & \le\frac{\lambda}{\sqrt{u-b_0(\gamma_n)}} & \text{ on } %
(b_0(\gamma_n), b_0(\gamma_n)+\varepsilon_0),\\
\intertext{and similarly,}
f_{\gamma_n}(u) & \le\frac{\lambda}{\sqrt{b_1(\gamma_n)-u}} & \text{ on } %
(b_1(\gamma_n)-\varepsilon_0, b_1(\gamma_n))
\end{align*}
To show continuity at $\gamma_0$, we now prove that $|\Phi_M(\gamma_n)-\Phi_M(\gamma_0)|\to 0$ as $n\to\infty$.
Unfortunately, to do this, we must consider separate cases.

Consider first the case where $b_0(\gamma_n)\nearrow b_0(\gamma_0)$. Define the positive numbers
$\delta_n=b_0(\gamma_0)-b_0(\gamma_n)+\frac14\varepsilon_0^2$ and
$\mu_n=b_1(\gamma_n)-b_1(\gamma_0)+\frac14\varepsilon_0^2$. Then
\begin{subequations}\label{PeriodEqn}
\begin{align}
\frac12|\Phi_M(\gamma_n)-\Phi_M(\gamma_0)| & \le \notag\\
\label{PeriodEqn1}&  \left|\int_{b_0(\gamma_n)}^{b_0(\gamma_n)+\delta_n}f_{\gamma_n}-%
\int_{b_0(\gamma_0)}^{b_0(\gamma_0)+\delta_n}f_{\gamma_0}\right|+\\
\label{PeriodEqn2} & \left|\int_{b_1(\gamma_n)-\mu_n}^{b_1(\gamma_n)}f_{\gamma_n}-%
\int_{b_1(\gamma_0)-\mu_n}^{b_1(\gamma_0)}f_{\gamma_0}\right|+\\
\label{PeriodEqn3} & \left|\int_{b_0(\gamma_n)+\delta_n}^{b_1(\gamma_n)-\mu_n}f_{\gamma_n}-%
\int_{b_0(\gamma_0)+\delta_n}^{b_1(\gamma_0)-\mu_n}f_{\gamma_0}\right|
\end{align}
\end{subequations}
We now turn our attention to each of the three terms in equations~\eqref{PeriodEqn}. For equation~\eqref{PeriodEqn1}, we have
\begin{gather}
\left|\int_{b_0(\gamma_n)}^{b_0(\gamma_n)+\delta_n}f_{\gamma_n}-%
\int_{b_0(\gamma_0)}^{b_0(\gamma_0)+\delta_n}f_{\gamma_0}\right| \notag\\
\label{PeriodEqn1A} \le\left|\int_{b_0(\gamma_n)}^{b_0(\gamma_0)}f_{\gamma_n}\right|+
\left|\int_{b_0(\gamma_0)}^{b_0(\gamma_n)+\delta_n}(f_{\gamma_n}-f_{\gamma_0})\right|+
\left|\int_{b_0(\gamma_n)+\delta_n}^{b_0(\gamma_0)+\delta_n}f_{\gamma_0}\right|
\end{gather}

We show that each of the terms in \eqref{PeriodEqn1A} can be made arbitrarily small. For the first term of
\eqref{PeriodEqn1A}:
$$
\int_{b_0(\gamma_n)}^{b_0(\gamma_0)}f_{\gamma_n}\le\int_{b_0(\gamma_n)}^{b_0(\gamma_0)}
\lambda\left[u-b_0(\gamma_n)\right]^{-1/2}=2\lambda(b_0(\gamma_0)-b_0(\gamma_n))^{1/2}<\lambda\varepsilon_0
$$
For the third term of \eqref{PeriodEqn1A}:
$$
\int_{b_0(\gamma_n)+\delta_n}^{b_0(\gamma_0)+\delta_n}f_{\gamma_0}\le%
\int_{b_0(\gamma_0)+\frac14\varepsilon_0^2}^{b_0(\gamma_0)+\delta_n}
\lambda\left[u-b_0(\gamma_0)\right]^{-1/2}=2\lambda[(\delta_n)^{1/2}-\varepsilon_0/2]
$$
which goes to $0$ as $n\to\infty$.
We handle the second term of \eqref{PeriodEqn1A} by applying the dominated convergence theorem: Note that
$b_0(\gamma_n)+\delta_n=b_0(\gamma_0)+\frac14\varepsilon_0^2$ and on the interval
$(b_0(\gamma_0),b_0(\gamma_0)+\frac14\varepsilon_0^2)$, $f_{\gamma_n}\to f_{\gamma_0}$ pointwise. 
Furthermore, on this interval 
$f_{\gamma_n}\le\lambda\left[u-b_0(\gamma_n)\right]^{-1/2}$ $<\lambda\left[u-b_0(\gamma_0)\right]^{-1/2}=g$, and
$\displaystyle\int_{b_0(\gamma_0)}^{b_0(\gamma_0)+\frac14\varepsilon_0^2} g=\lambda\varepsilon_0$. Thus, by dominated convergence
$$
\left|\int_{b_0(\gamma_0)}^{b_0(\gamma_n)+\delta_n}(f_{\gamma_n}-f_{\gamma_0})\right|\to 0\text{ as } n\to\infty
$$
For \eqref{PeriodEqn2} of equations \eqref{PeriodEqn} we write:
\begin{gather}
\left|\int_{b_1(\gamma_n)-\mu_n}^{b_1(\gamma_n)}f_{\gamma_n}-%
\int_{b_1(\gamma_0)-\mu_n}^{b_1(\gamma_0)}f_{\gamma_0}\right| \notag\\
\label{PeriodEqn1B} \le\left|\int_{b_1(\gamma_0)}^{b_1(\gamma_n)}f_{\gamma_n}\right|+
\left|\int_{b_1(\gamma_n)-\mu_n}^{b_1(\gamma_0)}(f_{\gamma_n}-f_{\gamma_0})\right|+
\left|\int_{b_1(\gamma_0)-\mu_n}^{b_1(\gamma_n)-\mu_n}f_{\gamma_0}\right|
\end{gather}
Arguing similarly, we conclude that each term of \eqref{PeriodEqn1B} can be made arbitrarily small. We omit the details.

Finally, for \eqref{PeriodEqn3} we have:
\begin{gather}
\left|\int_{b_0(\gamma_n)+\delta_n}^{b_1(\gamma_n)-\mu_n}f_{\gamma_n}-%
\int_{b_0(\gamma_0)+\delta_n}^{b_1(\gamma_0)-\mu_n}f_{\gamma_0}\right|\notag\\
\label{PeriodEqn3A} \le\left|\int_{b_0(\gamma_n)+\delta_n}^{b_0(\gamma_0)+\delta_n}f_{\gamma_n}\right|+
\left|\int_{b_0(\gamma_0)+\delta_n}^{b_1(\gamma_0)-\mu_n}(f_{\gamma_n}-f_{\gamma_0})\right|+
\left|\int_{b_1(\gamma_0)-\mu_n}^{b_1(\gamma_n)-\mu_n}f_{\gamma_n}\right|
\end{gather}
We now show that each of the terms in \eqref{PeriodEqn3A} can be made arbitrarily small. Note that
$b_0(\gamma_0)+\delta_n=b_0(\gamma_0)+[b_0(\gamma_0)-b_0(\gamma_n)]+\frac14\varepsilon_0^2<%
b_0(\gamma_0)+\frac12\varepsilon_0^2=b_0(\gamma_n)+[b_0(\gamma_0)-b_0(\gamma_n)]+\frac12\varepsilon_0^2<%
b_0(\gamma_n)+\frac34\varepsilon_0^2<b_0(\gamma_n)+\varepsilon_0$. Thus, for the first term of 
\eqref{PeriodEqn3A}:
\begin{align*}
\int_{b_0(\gamma_n)+\delta_n}^{b_0(\gamma_0)+\delta_n}f_{\gamma_n} \le &
\int_{b_0(\gamma_n)+\delta_n}^{b_0(\gamma_0)+\delta_n}
\lambda\left[u-b_0(\gamma_n)\right]^{-1/2}\\
= & 2\lambda\left[(2\delta_n-\varepsilon_0^2/4)^{1/2}-(\delta_n)^{1/2})\right]\to 0\text{ as } n\to\infty
\end{align*}
For the third term of \eqref{PeriodEqn3A}
we note that, similar to before, $b_1(\gamma_0)-\mu_n>b_1(\gamma_n)-\varepsilon_0$ thus:
\begin{align*}
\int_{b_1(\gamma_0)-\mu_n}^{b_1(\gamma_n)-\mu_n}f_{\gamma_n} \le &
\int_{b_1(\gamma_0)-\mu_n}^{b_1(\gamma_n)-\mu_n}
\lambda\left[b_1(\gamma_n)-u\right]^{-1/2}\\
= & -2\lambda\left[(\mu_n)^{1/2}-(2\mu_n-\varepsilon_0^2/4)^{1/2})\right]\to 0\text{ as } n\to\infty
\end{align*}
For the middle term of \eqref{PeriodEqn3A}, just note that $f_{\gamma_n}$ and $f_{\gamma_0}$ are both bounded on the
interval $(b_0(\gamma_0)+\delta_n,b_1(\gamma_0)-\mu_n)$ and that $f_{\gamma_n}\to  f_{\gamma_0}$ pointwise. Dominated convergence then implies that this term approaches zero as $n\to \infty$.

This is enough to verify continuity of the period function in the case when $b_0(\gamma_n)\nearrow b_0(\gamma_0)$.
We now complete the continuity proof by treating the case where $b_0(\gamma_n)\searrow b_0(\gamma_0)$. The proof here is essentially obtained by interchanging the roles of $\gamma_n$ and $\gamma_0$ in what has gone before. However, there are
some minor technical differences, which we point out.

To this end, define the positive numbers
$\delta_n=b_0(\gamma_n)-b_0(\gamma_0)+\frac14\varepsilon_0^2$ and
$\mu_n=b_1(\gamma_0)-b_1(\gamma_n)+\frac14\varepsilon_0^2$. Then
\begin{subequations}\label{PeriodEqn_SE}
\begin{align}
\frac12|\Phi_M(\gamma_0)-\Phi_M(\gamma_n)| & \le \notag\\
\label{PeriodEqn1_New}&  \left|\int_{b_0(\gamma_0)}^{b_0(\gamma_0)+\delta_n}f_{\gamma_0}-%
\int_{b_0(\gamma_n)}^{b_0(\gamma_n)+\delta_n}f_{\gamma_n}\right|+\\
\label{PeriodEqn2_New} & \left|\int_{b_1(\gamma_0)-\mu_n}^{b_1(\gamma_0)}f_{\gamma_0}-%
\int_{b_1(\gamma_n)-\mu_n}^{b_1(\gamma_n)}f_{\gamma_n}\right|+\\
\label{PeriodEqn3_New} & \left|\int_{b_0(\gamma_0)+\delta_n}^{b_1(\gamma_0)-\mu_n}f_{\gamma_0}-%
\int_{b_0(\gamma_n)+\delta_n}^{b_1(\gamma_n)-\mu_n}f_{\gamma_n}\right|
\end{align}
\end{subequations}

For equation~\eqref{PeriodEqn1_New}, we have
\begin{gather}
\left|\int_{b_0(\gamma_0)}^{b_0(\gamma_0)+\delta_n}f_{\gamma_0}-%
\int_{b_0(\gamma_n)}^{b_0(\gamma_n)+\delta_n}f_{\gamma_n}\right| \notag\\
\label{PeriodEqn1A_New} \le\left|\int_{b_0(\gamma_0)}^{b_0(\gamma_n)}f_{\gamma_0}\right|+
\left|\int_{b_0(\gamma_n)}^{b_0(\gamma_0)+\delta_n}(f_{\gamma_0}-f_{\gamma_n})\right|+
\left|\int_{b_0(\gamma_0)+\delta_n}^{b_0(\gamma_n)+\delta_n}f_{\gamma_n}\right|
\end{gather}

As before, we show that each of the terms in \eqref{PeriodEqn1A_New} can be made arbitrarily small. For the first term of
\eqref{PeriodEqn1A_New}:
$$
\int_{b_0(\gamma_0)}^{b_0(\gamma_n)}f_{\gamma_0}\le\int_{b_0(\gamma_0)}^{b_0(\gamma_n)}
\lambda\left[u-b_0(\gamma_0)\right]^{-1/2}=2\lambda(b_0(\gamma_n)-b_0(\gamma_0))^{1/2}<\lambda\varepsilon_0
$$
For the third term of \eqref{PeriodEqn1A_New}:
$$
\int_{b_0(\gamma_0)+\delta_n}^{b_0(\gamma_n)+\delta_n}f_{\gamma_n}\le%
\int_{b_0(\gamma_n)+\frac14\varepsilon_0^2}^{b_0(\gamma_n)+\delta_n}
\lambda\left[u-b_0(\gamma_n)\right]^{-1/2}=2\lambda[(\delta_n)^{1/2}-\varepsilon_0/2]
$$
which goes to $0$ as $n\to\infty$.
We now handle the second term of \eqref{PeriodEqn1A_New}. Since the functions
$f_{\gamma_n}$ are not defined on the entire domain of $f_{\gamma_0}$, there is a minor technical difference between this situation and the analogous one for \eqref{PeriodEqn1A}.
Define
$$ \hat{f}_{\gamma_n}=
\begin{cases} f_{\gamma_n} & \text{on $\left(b_0(\gamma_n),b_0(\gamma_n)+\frac14\varepsilon_0^2\right)$},\\
f_{\gamma_0} & \text{on $\big(b_0(\gamma_0),b_0(\gamma_n)\big)\cup \left(b_0(\gamma_n)+\frac14\varepsilon_0^2,b_0(\gamma_0)+\frac14\varepsilon_0^2\right)$}
\end{cases}
$$
Then $\displaystyle\int_{b_0(\gamma_n)}^{b_0(\gamma_0)+\delta_n}(f_{\gamma_0}-f_{\gamma_n})=%
\int_{b_0(\gamma_0)}^{b_0(\gamma_0)+\frac14\varepsilon_0^2}(f_{\gamma_0}-\hat{f}_{\gamma_n})$.

Then $\hat{f}_{\gamma_n}\to f_{\gamma_0}$ a.e. on $(b_0(\gamma_0),b_0(\gamma_0)+\frac14\varepsilon_0^2)$ and
$\hat{f}_{\gamma_n}\le g_n$ where
$$g_n=
\begin{cases} \lambda[u-b_0(\gamma_n)]^{-1/2} & \text{on $\left(b_0(\gamma_n),b_0(\gamma_n)+\frac14\varepsilon_0^2\right)$},\\
\lambda[u-b_0(\gamma_0)]^{-1/2} & %
\text{on $\big(b_0(\gamma_0),b_0(\gamma_n)\big)\cup \left(b_0(\gamma_n)+\frac14\varepsilon_0^2,b_0(\gamma_0)+\frac14\varepsilon_0^2\right)$}
\end{cases}
$$
Furthermore, on $(b_0(\gamma_0),b_0(\gamma_0)+\frac14\varepsilon_0^2)$, $g_n\to g$ a.e. where $g=\lambda[u-b_0(\gamma_0)]^{-1/2}$ and
$\displaystyle\lim_{n\to\infty}\int_{b_0(\gamma)}^{b_0(\gamma_0)+\frac14\varepsilon_0^2} g_n=\int_{b_0(\gamma)}^{b_0(\gamma_0)+\frac14\varepsilon_0^2} g$.
By a modified dominated convergence theorem which may be found, for example, in \cite{MR1013117}, we may conclude that
$$
\left|\int_{b_0(\gamma_n)}^{b_0(\gamma_0)+\delta_n}(f_{\gamma_0}-f_{\gamma_n})\right|=%
\left|\int_{b_0(\gamma_0)}^{b_0(\gamma_0)+\frac14\varepsilon_0^2}(f_{\gamma_0}-\hat{f}_{\gamma_n})\right|\to 0 \text{ as }n\to\infty
$$
Arguing similarly, we conclude that \eqref{PeriodEqn2_New} can be made arbitrarily small. We omit the details.
Lastly, for \eqref{PeriodEqn3_New} we have:
\begin{gather}
\left|\int_{b_0(\gamma_0)+\delta_n}^{b_1(\gamma_0)-\mu_n}f_{\gamma_0}-%
\int_{b_0(\gamma_n)+\delta_n}^{b_1(\gamma_n)-\mu_n}f_{\gamma_n}\right|\notag\\
\label{PeriodEqn3A_New} \le\left|\int_{b_0(\gamma_0)+\delta_n}^{b_0(\gamma_n)+\delta_n}f_{\gamma_0}\right|+
\left|\int_{b_0(\gamma_n)+\delta_n}^{b_1(\gamma_n)-\mu_n}(f_{\gamma_0}-f_{\gamma_n})\right|+
\left|\int_{b_1(\gamma_n)-\mu_n}^{b_1(\gamma_0)-\mu_n}f_{\gamma_0}\right|
\end{gather}
Arguing as we did for expression \eqref{PeriodEqn3A}, we may conclude that the first and third terms of \eqref{PeriodEqn3A_New} approach $0$ as $n\to\infty$.
For the middle term of \eqref{PeriodEqn3A_New}, proceed in the same way as we did to handle the second term of \eqref{PeriodEqn1A_New} by defining 
$$
\hat{f}_{\gamma_n}=
\begin{cases} f_{\gamma_n} & \text{on $\big(b_0(\gamma_n)+\delta_n,b_1(\gamma_n)-\mu_n\big)$},\\
f_{\gamma_0} & \text{on $\big(b_0(\gamma_0)+\delta_n,b_0(\gamma_n)+\delta_n\big)\cup \big(b_1(\gamma_n)-\mu_n,b_1(\gamma_0)-\mu_n\big)$}
\end{cases}
$$
and applying dominated convergence. This completes the proof of the continuity of the period function.
\end{proof}

\bibliography{ref}
\bibliographystyle{amsplain}

\end{document}